\documentclass[11pt,A4]{article} 
\usepackage{kimballnt}
\newcommand{\Ram}{\mathrm{Ram}}
\title{Central $L$-values and periods for $\GL(2)$}
\author{Kimball Martin}

\begin{document}

\maketitle

These are notes from 
{\em Automorphic Represenations, Automorphic Forms, L-functions and Related Topics, Jan. 21-25, 2008, RIMS, Kyoto, Conference Proceedings} with updated references.

In this talk we will attempt to give some sort of (biased) survey of formulas relating central values of twisted $\GL(2)$ $L$-functions to compact toric period integrals, ending with a summary of joint work with David Whitehouse.  First as a bit of motivation, and to try and indicate what is meant by the term ``period integral,'' let us recall the following two well-known examples.

In 1736, for $k \geq 1$, Euler proved the amazing formula
\begin{equation} \label{zeta2k}
\zeta(2k) = \sum_{n=1}^\infty \f 1{n^{2k}} = \f{2^{2k-1}\pi^{2k}|B_{2k}|}{(2k)!}
= \f{2^{2k-1}|B_{2k}|}{(2k)!} \left| \int_{-1}^1 \f{dx}{\sqrt{1-x^2}} \right|^{2k},
\end{equation}
where $B_k$ is the $k$-th Bernoulli number
\[ B_k = - \sum_{r=0}^{k-1} B_r \f 1{r!} \f{n!}{(n-k+1)!}. \]
The Birch--Swinnerton-Dyer conjecture predicts that if $E$ is an elliptic curve of rank 0 over $\Q$, then 
\begin{equation} \label{bsd}
 L(E,1) = \f{
\# \Sha(E)}
{(\#E_\tors(\Q))^2} \(\prod_p c_p(E) \)
\int_{E(\R)} |\omega| 
\end{equation}
where $\omega$ is a 
certain invariant (Neron) differential, 
$\Sha(E)$ denotes the Tate-Shafarevich group, and $c_p(E)$ are Tamagawa factors 
associated to the finitely many primes of bad reduction.  Both the period 
$\int |\omega|$ and the Tamagawa
factors $c_p$ are easy to compute.  On the other hand $\Sha(E)$ is not even
known to be finite.
We stress that the right-hand side of the equation
involves important arithmetic/geometric quantities.  (More generally, the same 
expression on 
the right is predicted for $\res_{s=1} L(E,s)/(2^r \text{Reg}(E))$ when the rank $r \neq 0$, where Reg$(E)$ denotes the regulator of $E$.)

The first formula is not of course a central value formula, but it is a special value formula.  
Here we wrote $\pi$ as
as the {\em period integral} $\int_{-1}^1 (1-x^2)^{1/2} dx$ so that we may realize (\ref{zeta2k}) as a formula relating the special values $\zeta(2k)$ with the integral of  a differential form on the open disc
of radius one in the plane over the {\em period} or closed geodesic $[-1,1]$.  For the formula (\ref{bsd}),
which is a central value formula,
the integral is taken over $E(\R)$, which will either be one or two closed cycles on the torus $E(\C)$, so
the integral on the right may also be called a {\em period integral}.  Classically the term period integral refers to an integral of a differential form over a (finite union of) closed geodesics or cycles.  One may extend this notion to higher dimensions by considering integrals over closed subgroups of algebraic groups.  

A general philosophy of how special values of $L$-functions should be related to period integrals
is given, up to rationality, by Deligne's conjectures.  More refined conjectures were later formulated by
Bloch-Kato, Beilinson and Ichino-Ikeda (see Ikeda's note in this volume for the latter).  We make no attempt say anything about this general philosophy, but focus on describing what has been done in the case of twisted $\GL(2)$ $L$-functions with a few words about the techniques involved.

\section{The setup}

We fix a number field $F$ and a cuspidal representation $\pi$ of $\GL_2(\A_F)$.  We will assume that
$\pi$ has trivial central character.
Now let $E$ be a quadratic extension of $F$ and $\chi$ be a unitary
character of $E^\times \A_F^\times \bs \A_E^\times$.  Our goal is to study the central value of the 
$L$-function
\[ L(s, \pi \times \chi):= L(s, \pi \times \pi(\chi)) = L(s, \pi_E \otimes \chi), \]
where $\pi(\chi)$ denotes the automorphic induction of $\chi$ from $\GL_1(\A_E)$ to a representation
of $\GL_2(\A_F)$ and $\pi_E$ denotes the base change of $\pi$ to $\GL_2(\A_E)$.  

In practice, one is
typically interested in varying the character $\chi$ and the extension $E$ and studying how the central values change.  For instance, note that if $\chi = 1$,
\[ L(s, \pi \times \chi) = L(s,\pi \times I_E^F(1)) = L(s,\pi)L(s,\pi \otimes \eta_{E/F}), \]
where $\eta_{E/F}$ denotes the quadratic character of $\A_F^\times$ corresponding by class field theory
to the non-trivial character of $\G(E/F)$.  Hence formulas for $L(1/2,\pi \times \chi)$ will yield
formulas for ratios 
\[ \f{L(1/2, \pi \otimes \eta_1)}{L(1/2,\pi \otimes \eta_2)}\] where $\eta_1$ and $\eta_2$ are quadratic characters of $\A_F^\times$ (see \cite{wald}---in fact, this was Waldspurger's original motivation for 
considering this situation).

So we will aim for a formula of the form
\[ L(1/2, \pi \times \chi) = c(\pi, E, \chi) |P(E,\chi)|^2 \]
where $c(\pi, E, \chi)$ is an ``arithmetic'' constant and $P(E,\chi)$ denotes some period integral.  We remark that one has a period squared and not just a period as the $L$-function is degree 4 over $F$ (or degree 2 over $E$).  Perhaps the first question one should ask is, 
what is an appropriate period to choose for $P(E,\chi)$?  In part this may depend on the application desired, but we will answer this question in a way which seems most aesthetically pleasing (and perhaps for this reason this answer will be well-suited for several applications).

We will make two assumptions now.  First, assume that $\epsilon(1/2, \pi \times \chi) = 1$.  The epsilon
factor is either $1$ or $-1$, and in the case it is -1, one automatically gets $L(1/2, \pi \times \chi) = 0$.  
(However if $\epsilon(1/2,\pi \times \chi) = -1$ one expects a similar formula in terms of a period for
a derivative of this $L$-function, which is a much harder question, cf.~\cite{GZ}.)  Secondly, we will assume
that $\pi$ and $\chi$ have disjoint ramification.  While this assumption is not needed to obtain a formula
for $L(1/2,\pi \times \chi)$, it will be required for our choice of {\em test vectors} for the period $P(E,\chi)$. 

Let $X$ be the set of pairs $(D, \pi_D)$ where $D$ is a quaternion algebra (possibly split, i.e., $M_2$) over $F$ containing $E$ such that $\pi_D$ corresponds to $\pi$ via the Jacquet-Langlands correspondence.
In particular $X$ always contains the element $(M_2, \pi)$ so $\pi$ need not actually come from a 
division algebra.  For $(D, \pi_D) \in X$, we choose $T$ to
be maximal torus $D$ such that $T(F) \iso E^\times$.  Then we define a period
\[ P_D : \pi_D \to \C \]
by
\[ P_D(\phi) = \int_{Z(\A_E)T(\A_F) \bs T(\A_E)} \phi(t) \chi\inv(t) dt. \]
To make a slightly more concrete connection with the geometry of the period integrals mentioned before,
we remark that $G$ modulo a compact open subgroup will be a hyperbolic 2- or 3- fold, and under this
image $T$ will be a finite collection of geodesics on the quotient.

\begin{thm} \rm{(\cite{wald}, \cite{jac-waldii})} \label{nonv-thm} $L(1/2, \pi \times \chi) \neq 0$ if and only if there exists a
$(D,\pi_D) \in X$ such that $P_D \neq 0$.
\end{thm}

\begin{thm}\rm{(\cite{tunnell})} \label{tunn-thm}
  The global linear functional $P_D$ on $\pi_D$ into a product of
local linear functionals $P_D = \prod_v P_{D,v}$, and 
there is a unique $(D, \pi_D) \in X$, characterized by local $\epsilon$-factors,
 such that $P_{D,v} \neq 0$ for all $v$.
\end{thm}

From now on, we fix this $D$ such that $P_{D,v}$ is nonvanishing for all $v$.  
Now the work of Gross-Prasad (\cite{GP}) singles out a nice vector $\phi \in \pi_D$, which is unique up to scaling, such that
\[ L(1/2, \pi \times \chi) \neq 0 \iff P_D(\phi) \neq 0. \]
This $\phi$ (up to scaling) will be called the {\em Gross-Prasad test vector}, as it can be used to test for
the existence of a non-zero $\chi$-equivariant linear form on $T$, or equivalently the nonvanishing of
the central $L$-value.
The test vector $\phi$ being ``nice'' essentially means that at all finite places $\phi_v$ is of minimum 
possible level such that $P_{D,v}(\phi_v) \neq 0$.  In particular if $\chi_v$ is unramified, $\phi_v$ is
a newform.  The work of Gross-Prasad is actually only done at the finite places, but at the infinite places
we may do the analogous thing:  choose $\phi_v$ to be a
 vector of minimum weight which transforms under $T$ by $\chi$.
Now it seems natural to look for a formula of the form
\[ L(1/2, \pi \times \chi) = c(\pi,E,\chi) \f{|P_D(\phi)|^2}{(\phi,\phi)} \]
where
\[ (\phi_1, \phi_2) = \int_{Z(\A_F)D^\times(F) \bs D^\times(\A_F)} \phi_1(g)\bar{\phi_2(g)} \]
and $\phi$ is the Gross-Prasad test vector.  We include the $(\phi, \phi)$ factor in the formula simply to
make the right hand side invariant under scaling of $\phi$.  There are two basic approaches to establishing such a formula: the theta correspondence and the relative trace formula.  We will summarize the results obtained by considering each approach in turn. 

\section{Theta correspondence results}

\begin{thm}\label{waldform} \rm{(\cite{wald})} For any $\phi \in \pi_D$,
\[ L(1/2, \pi \times \chi) \prod_{v \in S} \alpha_v(\pi, \phi, E, \chi) = \f 1{\zeta(2)}\f{|P_D(\phi)|^2}{(\phi,\phi)},\]
where $S$ is a finite set of ``bad'' primes, and $\alpha_v$ is a well-defined constant.
\end{thm}

\jump
{\bf Remarks:} 1.  Here the $\alpha_v$'s are defined in terms of integrals and may be 0.  The definition
of the $\alpha_v$ factors is sufficient for some applications such as rationality results, but is not
explicit enough for (direct application to) others, such as non-negativity of $L(1/2, \pi \times \chi)$ or
equidistribution results.

2.  To get a more explicit formula, one needs to choose a specific test vector $\phi$, and the Gross-Prasad test vector seems both natural and convenient for certain applications (Iwasawa theory, $p$-adic
$L$-functions---see \cite{vatsal}).

\jump
After Waldspurger's work, there were several efforts to refine his formula into something more explicit.  
Building off the work in \cite{GZ}, Gross proved the following result.

\begin{thm} \rm{(\cite{gross})} \label{gross-thm}  Suppose $F= \Q$, $E/F$ is a quadratic imaginary extension of discriminant $\Delta_E$, $\pi$ is
associated to a cusp form $f \in S_2(\Gamma_0(N))$ and $N$ is an inert prime in $E$.  Then 
\[ L(1/2, \pi \times \chi) = \f{N\sqrt{|\Delta_E|}}{(\pi |\mu(E)|)^2}  \f{|P_D(\phi)|^2}{(\phi,\phi)}, \]
where $\phi$ is the Gross-Prasad test vector, $\chi$ is unramified,
 and $\mu(E)$ is the group of roots of unity in $E$.
\end{thm}

Subsequently, 
Zhang (\cite{zhang}) generalized Gross's formula to the case of $F$ totally real, $E/F$ 
imaginary quadratic, $\pi$ corresponding to a holomorphic Hilbert modular form $f$ of parallel weight 2
and level $N$ assuming the conductor of $\chi$ is prime to $N$ and $\Delta_E$.  However Zhang is forced to choose a test vector which is not Gross-Prasad where $\chi$ ramifies.  With an assumption on $\chi_\infty$, Xue (\cite{xue}, \cite{xue2}) generalized Zhang's result to the case where $f$ is a holomorphic cusp form of arbitrary weight.  

The case of real quadratic extensions $E$ of $F=\Q$ and unramified characters $\chi$ was treated
by Popa (\cite{popa}) under the condition the level $N$ is squarefree and prime to $\Delta_E$.  Here,
the test vector $\phi$ used is the Gross-Prasad test vector.

\section{Relative trace formula results}

In this section we will just summarize the results coming from the relative trace formula approach,
but in the next section we will give a brief outline of the approach and make a few remarks.

The relevant relative trace formulas for this question were established by Jacquet in
\cite{jac-waldi} for the case of $\chi = 1$ (where a simpler trace formula is possible) and
in \cite{jac-waldii} for general $\chi$.  In these papers the aim of the relative trace formulas
is used to prove a nonvanishing criterion of compact toric periods on the quaternion algebra in terms of 
split toric periods on $\GL(2)$.  Jacquet suggested that these trace formula may be used to obtain another proof of Waldspurger's formula (Theorem \ref{waldform} above).  However, perhaps because
one already had Waldspurger's result, no real attempt was made to carry this out for some time.

The first real motivation for further progress was a to get a result on positivity of $L$-values.  In particular,
the Generalized Riemann Hypothesis implies that $L(1/2, \pi \times \chi) \geq 0$, but as mentioned earlier, Waldspurger's formula is not explicit enough to conclude this directly.   
In \cite{Guo}, Guo used the simpler relative trace formula of \cite{jac-waldi} to conclude
that $L(1/2, \pi) \geq 0$.  Generalizing work of \cite{Guo} in the case of $\chi = 1$,
Jacquet and Chen used the relative trace formula of \cite{jac-waldii} to prove the following result.

\begin{thm} \rm{(\cite{JC})} \label{jc-thm} 
For  $f \in C_c^\infty(D^\times(\A_F))$we have 
 \[ L(1/2, \pi \times \chi) = \f 12 \left(\prod_{v \in S} c_v \tilde J_v(f_v)\right) L(1, \pi, Ad) J_{\pi_D}(f), \]
 where $S$ is a finite set of bad places, $c_v$ is an explicit constant defined in terms of $L$-values
 and $\epsilon$-factors, and $\tilde J_v$ and $J_{\pi_D}$ are certain spectral distributions.
\end{thm}

The key point for us is that $J_{\pi_D}$ is be a sum of products of period integrals, which we will
expound upon in the next section.
(There is actually a slight restriction in \cite{JC} that $\pi$ not be dihedral with respect to $E$ for technical
simplifications, but this is later dealt with in the appendix to \cite{MW}.)  Using this result, the authors deduce $L(1/2, \pi \times \chi) \geq 0$.  (The positivity of central values is, of course, also immediate 
from the explicit formulas such as Theorem \ref{gross-thm} above or Theorem \ref{mw-thm} below.)

In joint work with David Whitehouse, we use the formula of Jacquet and Chen to prove the following

\begin{thm} \rm{(\cite{MW})} \label{mw-thm} For $\phi$ the Gross-Prasad test vector, we have
\begin{multline*}
 L(1/2, \pi \times \chi) = 2 \sqrt{\f{c(\chi)|\Delta_E|}{|\Delta_F|}} \( \prod_{v \in \Ram(\pi)} e_v(E/F) L(1,\eta_{E_v/F_v}) \right)\inv \\
\cdot \prod_{v | \infty} C_v(E,\pi,\chi)  \cdot  \f{L^S(1,\pi,Ad)L_S(1/2,\pi \times \chi)}{L_{\Ram(\chi)}(1,\eta)^2} \f{|P_D(\phi)|^2}{(\phi,\phi)},
\end{multline*}
where $\Ram(\circ)$ denotes the set of finite places where an object $\circ$ is ramified, 
 $S$ is the subset of places $v$ of $\Ram(\pi)$ such that either the conductor of $\pi_v$ is at least 2
or $v \in \Ram(E)$,
$\Delta_K$
denotes the discriminant of a field $K$, $e(L/K)$ denote the ramification degree of an extension $L/K$,
$c(\chi)$ is the absolute norm of the conductor of $\chi$, $\eta$ is the quadratic character attached to $E/F$, and $C_v(E,\pi,\chi)$ is an explicitly defined archimedean constant.
\end{thm}

Regarding the archimedean constants, if $v$ is real and split in $E$, then
\[ C_v = 2^{-k_v} \]
if $\pi_v$ is a discrete series of weight $k_v$ and
\[ C_v = \( \f{\lambda_v}{8\pi^2} \)^{\epsilon_v} \]
if $\pi_v = \pi(\mu_v, \mu_v\inv)$ is principal series with Laplacian eigenvalue $\lambda_v$ and
$\epsilon_v$ satisfies $\mu_v\chi_v = |\cdot|^{r_v}{\mathrm sgn}^{\epsilon_v}$. 
When $v$ is real inert or complex, the expressions for $C_v$ somewhat more complicated (involving
binomial coefficients and combinatorial products of parameters for $\pi$ and $\chi$),
and for those we refer to \cite{MW}.

For applications of these formula, we refer to $\cite{popa}$, $\cite{MW}$ and the references therein.

\section{The relative trace formula approach}

Now we will give a quick introduction to the relative trace formula approach.  Actually the way we
outline the method here is slightly more direct method than what is done in \cite{MW}, in that it does not involve going through Theorem \ref{jc-thm} above, though they are essentially equivalent.  However
we hope this presentation will help make the ideas in relative trace formula approach more transparent.  

We would like to remark that this relative trace formula approach can be used in complete generality (arbitrary central character and arbitrary ramification data for any given test function $\phi$).  
Specifically, the method used in \cite{MW} can deal with arbitrary ramification data and test function 
$\phi$ but makes use of Theorem \ref{jc-thm} which is only worked out for the central character of $\pi$ 
being equal to either 1 or $\eta_{E/F}$.  Of course, one should be able to extend Theorem \ref{jc-thm}
to an arbitrary central character, but that was proved with a somewhat different aim in mind.  In any case,
we would like to stress that to extend Theorem \ref{mw-thm} to greater generality should just boil down
to relatively simple local calculations.  However, for simplicity we will as above assume
trivial central character.

As before, let $T$ denote a torus over $F$ such that $T(F) \iso E^\times$.
Let $G$ and $G_\epsilon$ respectively denote the algebraic groups $\PGL(2)$ and $PD_\epsilon^\times$ over $F$ for an inner form $D_\epsilon$ of $\GL(2)$ containing $T$ 
parameterized by $\epsilon \in N_{E/F}(E^\times) \bs F^\times$.
Let $f$ and $f_\epsilon$ denote functions lying in $C_c^\infty(G(\A_E))$ and $C_c^\infty(G_\epsilon(\A_F))$, respectively.  We will assume that these functions are factorizable and at almost all
places are given by the characteristic function of the standard maximal compact subgroup.
Then we form the kernel function
\[ K(x,y) = K_f(x,y) := \sum_{\gamma \in G(E)} f(x\inv \gamma y) \]
associated to $f$ and one has the analogous kernel $K_\epsilon(x,y) = K_{f_\epsilon}(x,y)$ for $f_\epsilon$.  Note that
we are working with a kernel over $E$ for the group $G$, but a kernel over $F$ for the inner form $G_\epsilon$.  
Decomposition of the spectrum gives an identity of the form
\begin{equation} \label{kerexpn}
 \sum_{\gamma \in G(E)} f(x\inv \gamma y) = K(x,y) = \sum_{\Pi \textrm{ cusp}} K_\Pi(x,y) + K_{nc}(x,y),
 \end{equation}
where the sum on the right is over cuspidal representations $\Pi$ of $G(\A_E)$ and 
$K_{nc}(x,y)$ denotes
the non-cuspidal (residual plus continuous) part of the kernel.  The sum on the left is called the
geometric expansion of the kernel and the one on the right is called the spectral expansion.  The usual
Selberg trace formula (over $E$) 
is obtained by integrating these two expressions for $K(x,y)$ over the diagonal
subgroup $\Delta G \subset G \times G$.

The first relative trace formula, which will be a trace formula for $G$ over $E$, comes from integrating
\[ \int_{A(E) \bs A(\A_E)} \int_{G(F) \bs G(\A_F)} K(a,h)\chi(a)\eta_{E/F}(\det(g))dgda, \]
where $A$ denotes the diagonal split torus in $G$. In fact, the outer integral does not converge and so
 one needs to regularize it, but let us ignore convergence issues for the purpose of exposition and 
 merely say that a suitable truncation needs to be used to make sense of things.  From the geometric
 expansion of the kernel, one gets the the geometric side of the relative trace formula
 \[ \sum_{\gamma \in A(E) \bs G(E) / G(F)} I_\gamma(f) \]
 where $I_\gamma(f)$ is the relative orbital integral
\begin{equation} \label{splitorbint}
  I_\gamma(f) = \int_{A(\A_E)} \int_{G_\gamma(\A_F)} f(a\inv\gamma g)\chi(a)
 \eta_{E/F}(\det(g)) dg da, 
 \end{equation}
 and $G_\gamma(\A_F) = (\gamma\inv A(\A_E) \gamma \cap G(\A_F)) \bs G(\A_F)$.
Now let us consider the spectral side.
Using a Whittaker model for  a cuspidal representation $\Pi$ of $G(\A_E)$, one gets
the following (Proposition 4 of \cite{JC}):
\begin{align} \label{splitspec} \nonumber
 J_\Pi(f) :&= \int_{A(E) \bs A(\A_E)} \int_{G(F) \bs G(\A_F)} K_\Pi(a,h)\chi(a)\eta_{E/F}(\det(g))dgda \\
 \nonumber &=
\left( \sum_\phi \left( \prod_{v \in S} J_v(\phi; f)\right) \bar{\int \phi(g) \eta_{E/F}(\det(g))dg } \right) \cdot  L^S(1/2, \Pi \otimes \chi) \\
&=  \prod_{v \in S} \tilde J_v(\phi; f) \f{L_S(1,\eta_{E/F}) L^S(1/2, \Pi \otimes \chi)}{L^S(1,\pi,Ad)},
\end{align}
where $\phi$ runs over an orthonormal basis for $\Pi$, $S$ is a finite set of bad primes, and $J_v(\phi;f)$
and $\tilde J_v(\phi; f)$ are local spectral distributions given in terms of integral over the split
torus $A$.

The second relative trace formula, which will be a trace formula for $G_\epsilon$ over $F$, comes from integrating
\[ \int_{T(F) \bs T(\A_F)} \int_{T(F) \bs T(\A_F)} K_\epsilon(s,t) \chi(s\inv t) ds dt \]
where $T$ is a torus in $G_\epsilon$ such that $T(F) \iso E^\times$.  Analagous to (\ref{kerexpn}),
one writes out the geometric and spectral expansions of $K_\epsilon$ and then integrates each
side to obtain a trace formula.  The geometric side of this trace formula will
be
\[ \sum_{\gamma_\epsilon \in T(F) \bs G_\epsilon(F) / T(F)} I_{\epsilon, \gamma_\epsilon}(f_\epsilon) \]
where
\[ I_{\epsilon,\gamma_\epsilon}(f_\epsilon) = \int_{T(\A_F)} \int_{T_{\gamma_\epsilon}(\A_F)} f(s\gamma_\epsilon t) \chi(st) dsdt \]
and $T_\gamma = \gamma\inv T \gamma \cap T$.
On the spectral side, we have, for a cuspidal representation $\sigma$ of $G_\epsilon(\A_F)$,
\begin{align*} \label{quatspec}
\nonumber J_{\epsilon, \sigma}(f_\epsilon) :&= \int_{T(F) \bs T(\A_F)} \int_{T(F) \bs T(\A_F)} K_\epsilon(s,t) \chi(s\inv t) ds dt \\
\nonumber &= \int \int \sum_\phi (\sigma(f) \phi)(s)\bar{\phi(t)} \chi(s\inv t) ds dt \\
\nonumber &= \sum_{\phi} \int (\sigma(f)\phi)(s) \chi\inv(s)ds \bar{\int \phi(t)\chi\inv (t) dt},
\end{align*}
i.e.,
\begin{equation} \label{quatspec}
 J_{\epsilon, \sigma}(f_\epsilon) = \sum_{\phi} P_{D_\epsilon}(\sigma(f)\phi) \bar{P_{D_\epsilon}(\phi)},
\end{equation}
where $\phi$ runs over an orthonormal basis for $\sigma$ and $P_D$ denotes the period defined
earlier.

As is by now standard in this business, one attempts to match the (relative) 
orbital integrals on the geometric side
in order to get an equality of the spectral sides.  Recall that in the Selberg trace formula, the orbital
integrals are parameterized by conjugacy classes and thus one tries to find a matching of conjugacy classes on your two groups $G$ and $G'$ so that the individual orbital integrals will match up for
certain ``matching'' test functions $f$ and $f'$.

In the case of the relative trace formula, the integrals on the geometric sides are 
parameterized by double
cosets, as we have done above.
There is a notion of matching of these double cosets, however it involves varying the inner form $G_\epsilon$.
The point is that one matches ``regular'' double cosets $\gamma \in A(E) \bs G(E) /G(F)$ with pairs
$(\epsilon, \gamma_\epsilon)$ as $\epsilon$ runs over $N_{E/F}(E^\times) \bs F^\times$ and 
$\gamma_\epsilon$ runs over double cosets $T(F) \bs G_\epsilon(F) / T(F)$.  (The regular double cosets
will be those for which the orbital integral (\ref{splitorbint}) converges without any truncation.)  Given
this notion of matching $\gamma \corr (\epsilon, \gamma_\epsilon)$, we will say $f \in C_c^\infty(G(\A_E))$ matches with a family $(f_\epsilon)_\epsilon$ if the orbital integrals 
\[ I_\gamma(f) = I_{\epsilon, \gamma_\epsilon} (f_\epsilon) \]
match for any $\gamma \corr (\epsilon, \gamma_\epsilon)$.  Implicit here is that any matching family
$(f_\epsilon)$ will satisfy $f_\epsilon = 0$ for almost all $\epsilon$.  Jacquet (\cite{jac-waldii}) 
proves that for any such $f$ there exists a matching family $(f_\epsilon)$.  (The converse is not true;
for a converse, one would also need to sum over relative trace formulas of the first type where
the inner integration domain runs over inner forms of $G$ which split over $E$.)  Furthermore,
he regularizes and establishes these relative trace formulas and, by analysis of the continuous spectrum
and a standard argument, deduces that for matching $f \corr (f_\epsilon)$ one has
\begin{equation*}
J_{\pi_E}(f) = \sum_{\epsilon} J_{\epsilon, \pi_\epsilon}(f_\epsilon), 
\end{equation*}
where the sum on the right runs over $(D_\epsilon, \pi_\epsilon) \in X$, i.e., $\pi_\epsilon$'s vary over
Jacquet-Langlands transfers of $\pi$ to $D_\epsilon$ (including the trivial split case where 
$\epsilon = 1$).   By (\ref{quatspec}) and Theorem \ref{tunn-thm}, we know that for $\epsilon_0$ such that
$D_{\epsilon_0} = D$, the quaternion algebra picked out by Theorem \ref{tunn-thm}, 
\begin{equation} \label{specid}
 J_{\pi_E}(f) = J_{\epsilon_0, \pi'} (f_{\epsilon_0}). 
 \end{equation}
Now one can choose specific $f_{\epsilon_0}$ 
and $f$ such that the orbital integrals for $f$ and $f_{\epsilon_0}$ 
match for $\gamma \corr (\epsilon_0, \gamma_{\epsilon_0})$.  
Then one can put $f_{\epsilon_0}$ in a family $(f_\epsilon)$ which matches $f$ so one will have
(\ref{specid}) as above.  The point is to choose $f_{\epsilon_0}$ such that $\pi'(f_{\epsilon_0}): 
\pi' \to \pi'$ is orthogonal projection onto $\< \phi \>$ where $\phi$ is the Gross-Prasad 
(or your favorite) test vector in $\pi'$.  In this case, by (\ref{splitspec}) and (\ref{quatspec}), the
identity (\ref{specid}) becomes
\[  \prod_{v \in S} \tilde J_v(\phi; f) \f{L_S(1,\eta_{E/F}) L^S(1/2, \Pi \otimes \chi)}{L^S(1,\pi,Ad)} = \f {|P_D(\phi)|^2}{(\phi,\phi)}. \]
Hence it suffices to compute the local spectral distributions $\tilde J_v(\phi; f)$ for your specific choice of
$f$, which should be fairly simple for finite places $v \in S$.  For the infinite places, one uses
Barnes' Lemma from classical analysis and again this is not difficult.  This yields the desired formula.


\begin{thebibliography}{MW07}

\bibitem[GP91]{GP}
Benedict~H. Gross and Dipendra Prasad.
\newblock Test vectors for linear forms.
\newblock {\em Math. Ann.}, 291(2):343--355, 1991.

\bibitem[Gro87]{gross}
Benedict~H. Gross.
\newblock Heights and the special values of {$L$}-series.
\newblock In {\em Number theory (Montreal, Que., 1985)}, volume~7 of {\em CMS
  Conf. Proc.}, pages 115--187. Amer. Math. Soc., Providence, RI, 1987.

\bibitem[Guo96]{Guo}
Jiandong Guo.
\newblock On the positivity of the central critical values of automorphic
  {$L$}-functions for {${\rm GL}(2)$}.
\newblock {\em Duke Math. J.}, 83(1):157--190, 1996.

\bibitem[GZ86]{GZ}
Benedict~H. Gross and Don~B. Zagier.
\newblock Heegner points and derivatives of {$L$}-series.
\newblock {\em Invent. Math.}, 84(2):225--320, 1986.

\bibitem[Jac86]{jac-waldi}
Herv\'e Jacquet.
\newblock Sur un r\'esultat de {W}aldspurger.
\newblock {\em Ann. Sci. \'Ec. Norm. Sup\'er., IV.}, 19(2):185--229, 1986.

\bibitem[Jac87]{jac-waldii}
Herv\'e Jacquet.
\newblock Sur un r\'esultat de {W}aldspurger. {II}.
\newblock {\em Compositio Math.}, 63(3):315--389, 1987.

\bibitem[JC01]{JC}
Herv{\'e} Jacquet and Nan Chen.
\newblock Positivity of quadratic base change {$L$}-functions.
\newblock {\em Bull. Soc. Math. France}, 129(1):33--90, 2001.

\bibitem[MW09]{MW}
Kimball Martin and David Whitehouse.
\newblock Central ${L}$-values and toric periods for $\mathrm{GL}(2)$.
\newblock {\em Int. Math. Res. Not.}, 2009(1):141--191, 2009.

\bibitem[Pop06]{popa}
Alexandru~A. Popa.
\newblock Central values of {R}ankin {$L$}-series over real quadratic fields.
\newblock {\em Compos. Math.}, 142(4):811--866, 2006.

\bibitem[Tun83]{tunnell}
Jerrold~B. Tunnell.
\newblock Local {$\epsilon $}-factors and characters of {${\rm GL}(2)$}.
\newblock {\em Amer. J. Math.}, 105(6):1277--1307, 1983.

\bibitem[Vat04]{vatsal}
Vinayak Vatsal.
\newblock Special value formulae for {R}ankin {$L$}-functions.
\newblock In {\em Heegner points and Rankin $L$-series}, volume~49 of {\em
  Math. Sci. Res. Inst. Publ.}, pages 165--190. Cambridge Univ. Press,
  Cambridge, 2004.

\bibitem[Wal85]{wald}
Jean-Loup Waldspurger.
\newblock Sur les valeurs de certaines fonctions {$L$} automorphes en leur
  centre de sym\'etrie.
\newblock {\em Compositio Math.}, 54(2):173--242, 1985.

\bibitem[Xue06]{xue}
Hui Xue.
\newblock Central values of {R}ankin {$L$}-functions.
\newblock {\em Int. Math. Res. Not.}, Art. ID 26150, 41 pp, 2006.

\bibitem[Xue07]{xue2}
Hui Xue.
\newblock Central values of {$L$}-functions over {CM} fields.
\newblock {\em J. Number Theory}, 122(2):342--378, 2007.

\bibitem[Zha01]{zhang}
Shou-Wu Zhang.
\newblock Gross-{Z}agier formula for {${\rm GL}\sb 2$}.
\newblock {\em Asian J. Math.}, 5(2):183--290, 2001.

\end{thebibliography}

\end{document}